# Stabilization of Multi Fractional Order Differential Equation with Delay Time and Feedback Control

Badiea S. Hassoun[1], *
Raheam A. Al-Saphoy[2]

[1,2]College of Education for Pure Sciences, Tikrit University, Iraq
*badieass88@gmail.com

Sameer Q. Hassan[3]
[3]College of Education, Mustansiriyah University, Baghdad, Iraq

**Abstract:** The purpose of this article is to introduce the original results which devoted with the nonlinear control system problems involves of nonlinear differential equations of fractional orders. Thus, this system is described with a mixed of ordinary derivatives in the first and second order that, are unstable before feedback gain. More precisely, we investigate and analysis the nonlinear control system in related to feedback gain matrix. In addition, we prove that the considered system is locally asymptotically stabilizable via certain conditions. Then, this work reinforce through some application examples that programmed for illustrating and showing the stabilizability of the current systems with high efficiency and accuracy.

## 1. Introduction

The field of control and systems is currently one of the most important topics that play a good role in simplifying some systems. Thus, the control is involved in non-linear systems and the interpretation of complex phenomena, which is of great benefit in modernizing human civilization day after day [1].

Fractional calculus contributes to many important aspects such as science, engineering and physical applications. We mention some of its applications with fractional optimal control problems (FOCPs) that are subject to dynamic constraints with the objective function problems, in, bioscience [2], economic [3], and so on .

The stability of a nonlinear Langevin system of Mittag-Leffler (ML)-type fractional derivative affected by time-varying delays and differential feedback control stability has been studied by Zhao in ref. [4]. Then, Li and *et al.* are studying of the global stability problem for feedback control systems of impulsive fractional differential equations on networks [5]. Another direction studied by Qasim and *et* for some classes with composition FOCPs as in [6]. The stabilization and destabilization of fractional oscillators via a delayed feedback control has been considered by Čermák *et al.* in [7].

The Mittag-Leffler stabilization of fractional-order nonlinear systems with unknown control coefficients is verified and examined by Wang in [8].

　　The main objective of this work is to study non-linear systems with multiple fractional orders between zero and one with an ordinary derivative for control systems. The systems that are unstable were examined, then a feedback describes gain matrix the presence of control. So, after that we investigate and demonstrate the local stabilizabiliy with complete accuracy for nonlinear systems.

　　This outline of paper is organized as follows: Section 2, present some basic preliminaries concept and some auxiliary definitions. In section 3, we obtain the rigorous new results for the multiplying (fractional-one order ordinary) differential nonlinear feedback control system with some applications. Finally, we provide the results that have discovered which focused on the stabilizability problem of non-linear feedback control systems.

## 2. Preliminaries

　　In this section, we will present some important definitions and characterizations which play a good role to achieve the stabilizability concept of the considered system.

**Definition 2.1 [9]:** The formula for the Mittag-Leffler function is

$$E_\alpha(z) = \sum_{r=0}^{\infty} \frac{z^r}{\Gamma(r\alpha+1)}, \text{ where } z \in C \text{ and } \alpha > 0.$$

And the Mittag-Leffler function with two parameters:

$$E_{\alpha,\beta}(Z) = \sum_{r=0}^{\infty} \frac{z^r}{\Gamma(r\alpha+\beta)}, \text{ where } z \in C \text{ and } \alpha, \beta > 0.$$

**Definition 2.2 [10]:** The Gamma function is defined by the integral formula

$$\Gamma(Z) = \int_0^\infty t^{z-1} e^{-t} dt \; z \in C, (Re(z) > 0)$$

With the property of Gamma function

$$\Gamma(z+1) = z\Gamma(z).$$

**Definition 2.3 [10]:** The Laplace transform of the function $f(t)$, where $0 \le t < \infty$ and denoted $L\{f(t)\}$ as

$$L(f(t)) = F(s) = \int_0^\infty e^{-st} f(t) dt, s \in R^+$$

Where $R^+$ is the domain of nonnegative real numbers.

**Definition 2.4 [11]:** The point $x^* \in R^n$ is an equilibrium point for the differential equation

$$\dot{x} = f(t,x), \text{ if } f(t,x^*) = 0 \text{ for all } t.$$

**Definition 2.5 [12]:** The equilibrium point $x^*$ is said to be asymptotically locally stable if every solution starting sufficiently close to $x^*$ converges to $x^*$ as $t \to \infty$; that is when

$x(t, x_0) \to x^*$ as $t \to \infty$ for $x_0$ sufficiently close to $x^*$.

**Definition 2.6 [12]:** The equilibrium point $x^*$ is said to be asymptotically globally stable, if every solution converges to $x^*$ as $t \to \infty$; that is

$x(t, x_0) \to x^*$ as $t \to \infty$

Regardless of the initial point $x_0$ (or regardless of whether or not $x_0$ is close to $x^*$).

**Definition 2.7 [13]:** Let $\alpha > 0$ be a real number and let
$n = [\alpha] + 1$ for $\alpha \notin N_0$.
If $f(t) \in AC^n[a,b]$, then the Caputo fractional derivatives $aCD_t^\alpha f(t)$ and $tCD_b^\alpha f(t)$ exist almost everywhere on $[a,b]$.

$$aCD_t^\alpha f(t) = \frac{1}{\Gamma(n-\alpha)} \int_a^t (t-\tau)^{n-\alpha-1} f^{(n)}(\tau) d\tau, \qquad (1)$$

$$tCD_b^\alpha f(t) = \frac{(-1)^n}{\Gamma(n-\alpha)} \int_t^b (\tau-t)^{n-\alpha-1} f^{(n)}(\tau) d\tau, \qquad (2)$$

In particular, when $0 < \alpha < 1$ and $f(t) \in AC[a,b]$, then

$$aCD_t^\alpha f(t) = \frac{1}{\Gamma(1-\alpha)} \int_a^t (t-\tau)^{-\alpha} f'(\tau) d\tau, \qquad (3)$$

And

$$tCD_b^\alpha f(t) = \frac{-1}{\Gamma(1-\alpha)} \int_t^b (\tau-t)^{-\alpha} f'(\tau) d\tau. \qquad (4)$$

**Definition 2.8 [13]:** The power function and the constant function of the Caputo's derivative, is:
i. $aCD_t^\alpha (t-a)^\beta = \{0, \quad for\ \beta \in N_0\ and\ \beta < [\alpha],\ \frac{\Gamma(\beta+1)}{\Gamma(\beta-\alpha+1)} (t-a)^{\beta-\alpha},\ for\ \beta \in N_0\ and\ \beta \geq [\alpha]$, (5)
ii. $aCD_t^\alpha c = 0$, where $c$ is a constant (6)

**Lemma 2.9 [14]:** Let $A$ be an unbounded liner generator of a $C_\circ$-semigroup $T(t)$ satisfies following conditions
1. $\|T(t)\|_{L(x)} \leq Me^{wt}$,

2. If $A + B$ together with
$$D(A + B) = D(A)$$
is the generator of a perturbed $C_\circ$-semigroup $S(t)$ on $X$. And $B$ is abounded linear operator on $X$, satisfying
$\|S(t)\|_{L(x)} \leq Me^{(w+M\|B\|)t}$.

**Lemma 2.10 [15]:** Suppose $\alpha > 0$, $u(t)$ is a continuous nonnegative, non-decreasing function defined on $0 \leq t < T$, $u(t) < N$ (constant). Let $a(t)$ is a nonnegative and locally integrable on $0 \leq t < T$, and $Z(t)$ is a nonnegative function integrable on $0 \leq t < T$ with

$$a(t) \leq Z(t) + u(t) \int_0^t (t-s)^{\alpha-1} a(s)ds, \forall s \in [0,T].$$

Then,

$$a(t) \leq Z(t) + \int_0^t \left[\sum_{k=1}^{\infty} \frac{(\Gamma(\alpha)u(t))^k}{\Gamma(\alpha k)}(t-s)^{\alpha k-1}Z(s)\right]ds$$

In addition, if $Z(t)$ is a nondecreasing function on $0 \leq t < T$, we get

$$a(t) \leq Z(t)E_\alpha(\Gamma(\alpha)u(t)t^\alpha)$$

**Lemma 2.11 [16]:** There are finite real constants $N_1 \geq 1, N_2 \geq 1$, and $N_3 \geq 1$ for the Mittag-Leffler function $E_{\alpha,\beta}(At^\alpha)$. Then the following properties verified

i. For any $0 < \alpha < 1$, there are finite real constants $N_1$, $N_2$ such that

$$E_{\alpha,1}(At^\alpha) \leq N_1\|e^{At}\|, E_{\alpha,\alpha}(At^\alpha) \leq N_2\|e^{At}\|$$

where $A \in R^{n \times n}$.

ii. For any $\alpha \geq 1$, and $\beta = 1, 2, \alpha$ there are finite real constant $N_3$ such that

$$E_{\alpha,\beta}(At^\alpha) \leq N_3\|e^{At\alpha}\|.$$

## 3. Main result

In this section suppose that the following system

$$\dot{x} = Ax + \underline{K}(t-\tau)g\left(t, x(t), aCD_t^{\alpha_1}x(t), aCD_t^{\alpha_2}x(t)\right) + K\dot{x} \qquad (7)$$

With initial condition

$$x_0 = x(0)$$

When $[x_1, x_2, \ldots, x_n]^T \in R^{n \times 1}$ is the state vector of the system, $A = [a_{ij}]_{n \times n}$ is the constant matrix, and $g(.): R^+ \times R^n \to R^n$ is a continuous nonlinear function, $\underline{K}(t-\tau): R^+ \to R^n$ is represented the delay time. Then the system (7) is of type multi-first order and fractional dynamical nonlinear control system, where $aCD_t^{\alpha_1}$ and $aCD_t^{\alpha_2}$ are the two of Caputo fractional derivative with $\alpha_1, \alpha_2 \in (0,1)$.
$u(t) = K\dot{x}$ is the feedback, where $K \in R^{n \times n}$.

**Theorem 3.1:** The fractional nonlinear feedback control system (7) is locally asymptotically stable, if it satisfies the following conditions:

i)　　$Re\ (eig\ (I-K)^{-1}A) < 0$, and
　　　$\omega = -Re(eig\ (I-K)^{-1}A) > M_3\ max\{(I-K)^{-1}\}$,
　　Where,
　　　$M_3 = MM_1M_2 \in R^+$.

ii)　　$\|\underline{K}(t-\tau)\| \leq M_1$, $M_1 \in R^+$.

iii)　　$\|g\left(t, x(t), aCD_t^{\alpha_1}x(t), aCD_t^{\alpha_2}x(t)\right)\| \leq M_2(\|x(t)\| + \|aCD_t^{\alpha_1}x(t)\| + \|aCD_t^{\alpha_2}x(t)\|), M_2 \in R^+$.

**Proof:**
From (7), we get

$$(I-K)\dot{x} = Ax + \underline{K}(t-\tau)g\left(t, x(t), aCD_t^{\alpha_1}x(t), aCD_t^{\alpha_2}x(t)\right),$$

$$\dot{x} = (I - K)^{-1}\left[Ax + \underline{K}(t - \tau)g\left(t, x(t), aCD_t^{\alpha_1}x(t), aCD_t^{\alpha_2}x(t)\right)\right],$$
(8)

Applying the Laplace transform in equation (8), we have
$$sX(s) - x_0 = (I - K)^{-1}AX(s) + (I - K)^{-1}$$
$$\cdot L\left\{\underline{K}(t - \tau)g\left(t, x(t), aCD_t^{\alpha_1}x(t), aCD_t^{\alpha_2}x(t)\right)\right\},$$

Therefore,
$$X(s) = (s - (I - K)^{-1}A)^{-1}x_0 + (s - (I - K)^{-1}A)^{-1}$$
$$\cdot \left((I - K)^{-1}L\left\{\underline{K}(t - \tau)g\left(t, x(t), aCD_t^{\alpha_1}x(t), aCD_t^{\alpha_2}x(t)\right)\right\}\right),$$

By using the Laplace inverse transform for both sides, we get
$$\|x(t)\| \leq \|e^{(I-K)^{-1}At}x_0\| + \left\|\int_0^t e^{(I-K)^{-1}A(t-\sigma)}(I - K)^{-1}\right.$$
$$\left. \cdot \underline{K}(\sigma - \tau)g\left(\sigma, x(\sigma), aCD_t^{\alpha_1}x(\sigma), aCD_t^{\alpha_2}x(\sigma)\right)d\sigma\right\|,$$

Now, by using lemma (2.9), we get
$$\|x(t)\| \leq Me^{-\omega t}\|x_0\| + (I - K)^{-1}\int_0^t Me^{-\omega(t-\sigma)}$$
$$\cdot \|\underline{K}(\sigma - \tau)\|\|g\left(\sigma, x(\sigma), aCD_t^{\alpha_1}x(\sigma), aCD_t^{\alpha_2}x(\sigma)\right)\|d\sigma,$$

By condition (ii), we have that
$$\|x(t)\| \leq Me^{-\omega t}\|x_0\| + (I - K)^{-1}\int_0^t Me^{-\omega(t-\sigma)}$$
$$\cdot M_1\|g\left(\sigma, x(\sigma), aCD_t^{\alpha_1}x(\sigma), aCD_t^{\alpha_2}x(\sigma)\right)\|d\sigma, M_1 \in R^+.$$

And by using condition (iii), we get

$$\|x(t)\| \leq Me^{-\omega t}\|x_0\| + (I - K)^{-1}M_3\int_0^t e^{-\omega(t-\sigma)}$$
$$\cdot \left(\|x(\sigma)\| + \|aCD_t^{\alpha_1}x(\sigma)\| + \|aCD_t^{\alpha_2}x(\sigma)\|\right)d\sigma,$$

Where,
$$M_3 = MM_1M_2 \in R^+.$$

Let
$$k_1(t) = \|aCD_t^{\alpha_1}x(t)\|$$

And
$$k_2(t) = \|aCD_t^{\alpha_2}x(t)\|$$

Where $k_1(t)$ and $k_2(t)$ can describe as follows:

$$k_2(t) = \left\|\frac{1}{\Gamma(1-\alpha_2)}\int_0^t (t - \mu)^{\alpha_2}\dot{x}(\mu)d\mu\right\|$$

and
$$\|x(t)\| \leq Me^{-\omega t}\|x_0\| + (I - K)^{-1}M_3\left(\int_0^t e^{-\omega(t-\sigma)}\|x(\sigma)\|d\sigma\right.$$
$$+ \int_0^t k_1(t)e^{-\omega(t-\sigma)}d\sigma + \int_0^t k_2(t)e^{-\omega(t-\sigma)}d\sigma\right),$$
$$\leq Me^{-\omega t}\|x_0\| + (I - K)^{-1}M_3\left(\int_0^t e^{-\omega(t-\sigma)}\|x(\sigma)\|d\sigma\right.$$
$$+ k_1(t)e^{-\omega t}\int_0^t e^{\sigma t}d\sigma + k_2(t)e^{-\omega t}\int_0^t e^{\sigma t}d\sigma\right),$$
$$\leq Me^{-\omega t}\|x_0\| + (I - K)^{-1}M_3\left(\int_0^t e^{-\omega(t-\sigma)}\|x(\sigma)\|d\sigma\right.$$
$$+ \frac{(1-e^{-\omega t})}{\omega}(k_1(t) + k_2(t))\right),$$
$$\|x(t)\| \leq Me^{-\omega t}\|x_0\| + \frac{M_3(1-e^{-\omega t})}{\omega}(I - K)^{-1}(k_1(t) + k_2(t))$$

$$+(I - K)^{-1}M_3 \int_0^t e^{-\omega(t-\sigma)}\|x(\sigma)\|d\sigma ,  \quad (9)$$

Multiplying the equation (9) by the term ($e^{\omega t}$), to get

$$\| x(t) \| e^{\omega t} \leq M\|x_0\| + \frac{M_3(e^{\omega t}-1)}{\omega}(I - K)^{-1}(k_1(t) + k_2(t))$$
$$+(I - K)^{-1}M_3 \int_0^t e^{\omega \sigma}\|x(\sigma)\|d\sigma ,$$

Now, by using lemma (2.10), we get

$$\| x(t) \| e^{\omega t} \leq (M\|x_0\| + \frac{M_3(e^{\omega t}-1)}{\omega}$$
$$\cdot (I - K)^{-1}(k_1(t) + k_2(t)))e^{M_3(I-K)^{-1}t} ,  \quad (10)$$

Therefore, by multiplying the equation (10) by the term ($e^{-\omega t}$), becomes

$$\| x(t) \| \leq (M\|x_0\| + \frac{M_3(e^{\omega t}-1)}{\omega}(I - K)^{-1}(k_1(t) + k_2(t)))e^{(M_3(I-K)^{-1}-\omega)t} ,$$

When
$$t \to \infty, \|x(t)\| \to 0$$
For $\omega > M_3 \, max\{(I - K)^{-1}\}$.
Consequently, the system (1) is asymptotically locally stable. ■

**Example 3.2:** Consider the following differential nonlinear without feedback control system:

$$\dot{x} = Ax + \underline{K}(t - \tau)g\left(t, x(t), aCD_t^{\alpha_1}x(t), aCD_t^{\alpha_2}x(t)\right),  \quad (11)$$

This nonlinear control system consists of nonlinear differential equations of fractional orders with a mixed of ordinary derivatives of the first order that are unstable before feedback gain matrix.

Thus, we examine this nonlinear control system after applying feedback gain matrix and prove the asymptotic local stabilizability of the system by using the conditions in the theorem (3.1).

$$\dot{x}_1 = -15x_1 + 15x_2 \quad (12)$$
$$\dot{x}_2 = 110x_1 - x_2 + x_2x_1^{\frac{4}{3}} \quad (13)$$
$$\dot{x}_3 = x_1x_2^{\frac{2}{3}} - 3x_3 \quad (14)$$

Where,
$A = [-15\ 15\ 0\ 110\ -1\ 0\ 0\ 0\ -1]$,
$\underline{K}(t - \tau) = t - 1$,
$x(t) = [x_1(t)\ x_2(t)\ x_3(t)]$,

And
$$g\left(t, x(t), aCD_t^{\alpha_1}x(t), aCD_t^{\alpha_2}x(t)\right) = \frac{1}{t-1}\left[0\ x_2x_1^{\frac{4}{3}}\ x_1x_2^{\frac{2}{3}}\right]$$

By using the definition (2.5) of power function
$$aCD_t^{\alpha}x^{\beta} = \frac{\Gamma(1+\beta)}{\Gamma(1+\beta-\alpha)}x^{\beta-\alpha}, \text{ we get}$$

i) When $\alpha_1 = \frac{2}{3}$ and $f(x) = x_1^2$

Then
$$aCD_t^{\frac{2}{3}}x_1^2 = \frac{\Gamma(3)}{\Gamma(3-\frac{2}{3})}x_1^{2-\frac{2}{3}} = \frac{2}{\Gamma(\frac{7}{3})}x_1^{\frac{4}{3}} = 1.68x_1^{\frac{4}{3}}$$

ii) When $\alpha_2 = \frac{3}{5}$ and $f(x) = x_2$

Then
$$aCD_t^{\frac{3}{5}}x_2 = \frac{\Gamma(2)}{\Gamma(2-\frac{3}{5})}x_2^{1-\frac{3}{5}} = \frac{1}{\Gamma(\frac{7}{5})}x_2^{\frac{2}{5}} = 1.127x_2^{\frac{2}{5}}$$

Thus,
$$\underline{K}(t-\tau)g\left(t,x(t),aCD_t^{\alpha_1}x(t),aCD_t^{\alpha_2}x(t)\right) = (t-1)\frac{1}{t-1}\left[0 \; x_2x_1^{\frac{4}{3}} \; x_1x_2^{\frac{2}{3}}\right]$$
$$= \left[0 \; x_2x_1^{\frac{4}{3}} \; x_1x_2^{\frac{2}{3}}\right]$$

Now, the figures 3.1: (a), (b) and (c) show the nonlinear fractional control system (11) is unstable before applying feedback gain matrix, as follows:

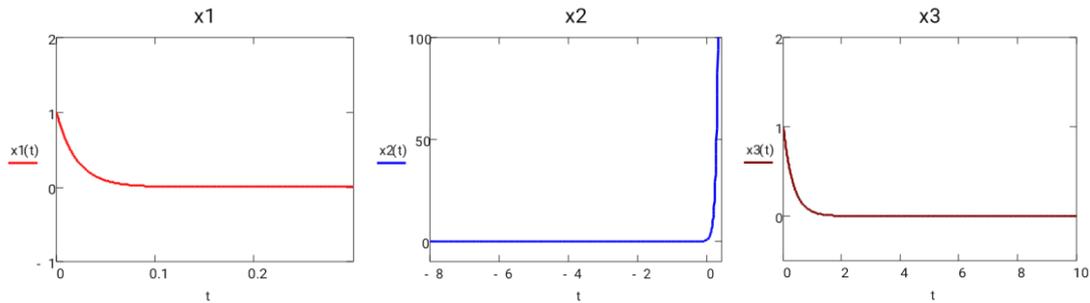

**Fig. 3.1**: Show the unstable nonlinear fractional control system.
a) The state solution $x_1$ of the nonlinear fractional control system (11) without feedback.
b) The state solution $x_2$ of the nonlinear fractional control system (11) without feedback.
c) The state solution $x_3$ of the nonlinear fractional control system (11) without feedback.

Now, we study the stabilizability notion of the nonlinear control system after applying feedback gain control input
$$u(t) = K\dot{x}$$
For system (1), as follows:
$$\dot{x} = (I-K)^{-1}\left[Ax + \underline{K}(t-\tau)g\left(t,x(t),aCD_t^{\alpha_1}x(t),aCD_t^{\alpha_2}x(t)\right)\right]$$
(15)
The feedback matrix was chosen, as follows:
$$K = [1 \; 15 \; 0 \; 110 \; 0 \; 0 \; 0 \; 0 \; -3\,],$$

$$I - K = [0 \; -15 \; 0 \; -110 \; 1 \; 0 \; 0 \; 0 \; 4\,],$$

$$(I-K)^{-1} = [-0.0006 \; -0.0009 \; 0 \; 1.0005 \; -1.0005 \; 0 \; 0 \; 0 \; -0.25\,],$$
And,
$$(I-K)^{-1}A = [-0.981 \; 0 \; 0 \; 1.0005 \; -1.0005 \; 0 \; 0 \; 0 \; -0.25\,],$$
then,
$$(I-K)^{-1}\underline{K}(t-\tau)g\left(t,x(t),aCD_t^{\alpha_1}x(t),aCD_t^{\alpha_2}x(t)\right)$$
$$= [-0.0006 \; -0.0009 \; 0 \; 1.0005 \; -1.0005 \; 0 \; 0 \; 0 \; -0.25\,]\left[0 \; x_2x_1^{\frac{4}{3}} \; x_1x_2^{\frac{2}{3}}\right] = \left[-0.0009x_2x_1^{\frac{4}{3}} \; -1.0005x_2x_1^{\frac{4}{3}} \; -0.25x_1x_2^{\frac{2}{3}}\right]$$

The nonlinear control system (25) with feedback takes the following form:

$$\dot{x}_1 = -0.981x_1 - 0.0009 x_2 x_1^{\frac{4}{3}} \tag{16}$$
$$\dot{x}_2 = 1.0005(x_1 - x_2) + x_2 x_1^{\frac{4}{3}} - 1.0005 x_2 x_1^{\frac{4}{3}} \tag{17}$$
$$\dot{x}_3 = x_1 x_2^{\frac{2}{3}} - 0.25 x_3 - 0.25 x_1 x_2^{\frac{2}{3}} \tag{18}$$

Then, prove and satisfy the asymptotic local stabilizability of the system (11) by using the conditions in the theorem (3.1), as follows:

i) $\quad Re\,(eig\,(I-K)^{-1}A) = -0.981, -1.0005, -0.25$

And
$$\omega = -Re(eig\,(I-K)^{-1}A) = 0.25 > M_3\,max\{(I-K)^{-1}\},$$
$$= \frac{1}{2}(-0.25) = -0.125$$

Where $M_3 = \frac{1}{2} \in R^+$.

ii) $\quad \| \underline{K}(t-\tau) \| = \| t-1 \|, \quad 0 < \tau < t$
$$\leq \| t \| \in R^+$$

iii) $\quad \| g\left(t, x(t), aCD_t^{\alpha_1}x(t), aCD_t^{\alpha_2}x(t)\right) \| = \| \frac{1}{t-1} \| \| \sqrt{(x_2 x_1^{\frac{4}{3}})^2 + (x_1 x_2^{\frac{2}{3}})^2} \|$
$$= \lim_{t \to \infty} \| \frac{1}{t-1} \| \| \sqrt{(x_2 x_1^{\frac{4}{3}})^2 + (x_1 x_2^{\frac{2}{3}})^2} \|$$
$$= 0 \text{ , when } t \to \infty \text{ and } M_2 = 1 \in R^+$$

$\| g\left(t, x(t), aCD_t^{\alpha_1}x(t), aCD_t^{\alpha_2}x(t)\right) \|$
$$\leq M_2(\|x(t)\| + \|aCD_t^{\alpha_1}x(t)\| + \|aCD_t^{\alpha_2}x(t)\|), M_2 \in R^+.$$

Thus, the figures 3.2: (a), (b) and (c) display the simulation result of theorem (3.1), which proves that the zero solution of the nonlinear system in (11) is asymptotically locally stabilizable.

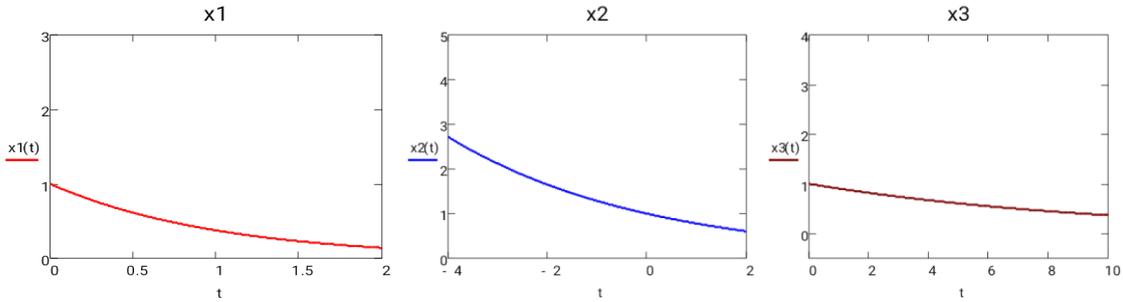

**Fig. 3.2:** Show the stabilizable nonlinear fractional control system.

a) The state solution $x_1$ of the nonlinear fractional control system (11) with feedback.
b) The state solution $x_2$ of the nonlinear fractional control system (11) with feedback.
c) The state solution $x_3$ of the nonlinear fractional control system (11) with feedback.

## 4. Conclusions

A new outcomes has been explored in nonlinear dynamical systems for double fractional with ordinary order in this paper related to some necessary conditions. Then, the stabilizability of nonlinear systems with control was obtained by feedback gain for the nonlinear systems class. So that, the precise results are obtained in locally case according to the conditions of theorem 3.1, which have demonstrated their accuracy in some applications. Also, as the work technology has been programmed and reinforced with the illustrative examples that have shown the efficiency of the stabilizability of the considered systems. Finally, may be interested to extend the obtained results in this work to the case of regional observer problem in distributed parameter systems as in [17-18].